\documentclass{amsart}
\usepackage{amssymb}
\usepackage{amsthm}
\usepackage{amsmath}
\usepackage{bm}
\usepackage{enumerate}
\usepackage[dvips]{graphics,color}

\theoremstyle{plain}
\newtheorem{Thm}{Theorem}[section]
\newtheorem{Lem}[Thm]{Lemma}

\newtheorem{Cor}[Thm]{Corollary}

\newcommand{\NN}{{\mathbb N}}

\newcommand{\RR}{{\mathbb R}}

\newcommand{\BB}{\mathbb B}

\newcommand\Span{{\mathop{\mbox{\rm span }}}}

\newcommand{\Inj}{\mathfrak{I}}
\newcommand{\Int}[1]{\stackrel{\circ}{#1}}

\def\<{\langle}
\def\>{\rangle}

\title{A note on James spaces and superstrictly singular operators}
\author{Isabelle Chalendar, Emmanuel Fricain and Dan Timotin}

\begin{document}

\begin{abstract}
An elementary lemma is used in order to show that the natural inclusion $J_p\to J_q$ of James spaces is superstrictly singular for $p<q$. As a consequence, it is shown that an operator without nontrivial invariant subspaces constructed by Charles Read is superstrictly singular. 
\end{abstract}

\address[I. Chalendar, E. Fricain]{Universit\'e de Lyon, Universit\'e Lyon 1, I. C. J., UFR de
Math\'ematiques, 43 bld. du 11 novembre 1918,
69622 Villeurbanne Cedex, France}

\email{chalenda@math.univ-lyon1.fr, fricain@math.univ-lyon1.fr}

\address[D. Timotin]{Institute of Mathematics of the 
Romanian Academy, PO Box 1-764, Bucharest 014700, Romania}

\email{Dan.Timotin@imar.ro}

\maketitle

\section{Introduction}

An operator $T:X\to Y$ between two Banach spaces is called
\begin{enumerate}
\item\emph{compact} if the image of the unit ball is relatively compact;
\item\emph{strictly singular} \cite{LT} if there is no infinite-dimensional subspace $E$ of $X$ such that the restriction $T|E$ is bounded below;

\item
\emph{superstrictly singular} \cite{M, P} if for every $\epsilon>0$ there exists $n\in\NN$ such that 
\[
\inf_{x\in E, \|x\|= 1}\|Tx\|<\epsilon
\]
for every $n$-dimensional subspace $E$ of $X$. 
\end{enumerate}
It is easy to see that any compact operator is superstrictly singular, and any superstrictly singular operator is strictly singular. Also, on a Hilbert space the three classes coincide. In general, they are closed ideals; also, they are preserved under formation of finite direct sums. These properties can be found, for instance, in~\cite{LT} for the first two classes and in~\cite{M} (see also~\cite{P}) for the third class. Note that superstrictly singular are sometimes called in the literature finitely strictly singular~\cite{SSTT}.

In~\cite{P, SSTT} it is shown that the natural inclusion of $\ell^p$ into $\ell^q$ ($p<q$) is superstrictly singular. This fact is used in~\cite{SSTT} for the investigation of the ideal structure of $L(\ell^p\oplus\ell^q)$.

On the other hand, Charles Read constructs in~\cite{R} an example of a strictly singular operator without nontrivial closed invariant subspaces (this seems to be an older question of Pe\l czy\'nski). Read's operator acts on an infinite direct sum of James spaces. The \emph{James space} $J_p$ ($0<p<\infty$) is introduced for $p=2$ in~\cite{J} and for other values of $p$ in~\cite{R}. It is the subspace of $c_0$ formed by those vectors $a=(a_i)_{i=1}^\infty$ such that 
\[
\|a\|_{J_p}=\sup\big\{ \big(\sum_{j=2}^n |a_{i_j}-a_{i_{j-1}}|^p \big)^{1/p} : i_1<\dots<i_n,\quad n\in\NN \big\}<\infty.
\]
A short account of the basic properties of $J_2$ can be found in~\cite{S}. Most notably (endowed with an equivalent norm) it has been  the first example of a non-reflexive Banach space that is  isometrically isomorphic to its second dual. 

The starting point of the construction in~\cite{R} is the observation that the natural inclusion of $J_p$ into $J_q$ ($p<q)$ is strictly singular. Our main result below (Theorem~\ref{th:main}) states that this natural inclusion is in fact even superstrictly singular. As an application, we show that the operator without nontrivial invariant subspaces constructed in~\cite{R} is actually superstrictly singular.

Our basic tool will be an interesting  elementary result (Lemma~\ref{le:pm1}), which might find applications also in other areas. Namely, it is not hard to prove that if $X\subset c_0$ is a $k$-dimensional subspace, then there is a vector $x\in X$ of norm 1, with $k$ coordinates equal in modulus to~1 (see, for instance,~\cite{M}). We show that one can actually ensure that these $k$ coordinates have alternating signs.

The authors thank Charles Read for bringing this problem to their attention.

\section{Main results}

We rely on the following technical result, whose proof will be given in section~\ref{se:proof}.

\begin{Lem}\label{le:pm1}
If $X\subset c_0$ is a subspace of dimension $k$, then there is $x\in X$, $\|x\|=1$ and indices $i_1 <\dots <i_k$, such that $x_{i_j}=(-1)^j$.
\end{Lem}

In order to show that for $p<q$ the natural inclusion from $J_p$ into $J_q$ is superstrictly singular, we will adapt the proof of Proposition~3.3 in~\cite{SSTT}. The main point is that we have to use Lemma~\ref{le:pm1} above instead of the simpler lemma from~\cite{M}.

\begin{Thm}\label{th:main}
If $p<q$, then the natural inclusion $\Inj_{p,q}:J_p \to J_q$ is superstrictly singular.
\end{Thm}

\begin{proof} %Note first that it is enough to prove the result in the case of real-valued James spaces (which is actually the case originally considered in ~\cite{J}).
%Indeed, the complex-valued James space can also be considered as a real James space, with an equivalent norm (up to the constant $\sqrt{2}$), and the formal identity remains the same.
For any $x\in J_p$ we have 
\[
|x_{i+1}-x_i|^q\le (2\|x\|_\infty)^{q-p}|x_{i+1}-x_i|^p,
\]
and therefore 
\[
\|\Inj_{p,q}x\|_{J_q}\le 2^{\frac{q-p}{p}} \|x\|_\infty^{\frac{q-p}{p}} \|x\|_{J_p}^\frac{p}{q}.
\]

Suppose that $E$ is a subspace of $J_p$ with $\dim E=k$. Since $J_p\subset c_0$, Lemma~\ref{le:pm1} yields the existence of $x\in E$, $\|x\|_\infty=1$ and indices $i_1<\dots< i_k$ such that 
$x_{i_j}=(-1)^j$. This implies that 
\[
\|x\|_{J_p}\ge \big((k-1)2^p \big)^{1/p}=2(k-1)^{1/p}.
\]
Then, if $z=\frac{x}{\|x\|_{J_p}}$, then $\|z\|_{J_p}=1$, $\|z\|_\infty\le \frac{1}{2(k-1)^{1/p}}$, whence
\[
\|\Inj_{p,q}z\|_{J_q}\le 2^{-p/q} \frac{1}{(k-1)^\frac{q-p}{pq}}.
\]
It is now clear that for any $\epsilon>0$ we can find $k$ such that the quantity in the right hand side is smaller than $\epsilon$. This proves that $\Inj_{p,q}$ is superstrictly singular.
\end{proof}

We apply now Theorem~\ref{th:main} to discuss an operator considered by Charles Read. 
In~\cite{R}, one constructs on the Banach space $X=\ell_2\oplus X_1$, where
\[
X_1=(\bigoplus_{i=1}^\infty J_{p_i})_{\ell^2}
\]
($p_i$ a strictly increasing sequence of integers), an operator $T$ 
which is strictly singular and has no nontrivial invariant subspaces. The construction is rather intricate: one finds a certain basis $(e_i)_{i=0}^\infty$ and define $T$ by the conditions $Te_i=e_{i+1}$; there is then a lot of work to show that $T$ extends to a continuous operator that is strictly singular and has no nontrivial invariant subspaces.

It turns out however that with the aid of Theorem~\ref{th:main} one can show that the operator $T$ is actually superstrictly singular. Indeed, in~\cite{R} one shows that $T$ is a compact perturbation of $(0\oplus W_1)$ , with $W_1:X_1\to X_1$ a weighted unilateral shift with  weights tending to 0; thus
\[
W_1((x_1, x_2,x_3,\dots))=(0, \beta_1 x_1, \beta_2 x_2,\dots)
\]
with $\beta_i\to 0$. Note that one should rather write $\beta_i \Inj_{p_i,p_{i+1}}x_i$ instead of $\beta_i x_i$.

\begin{Lem}
The operator $W_1$ is superstrictly singular.
\end{Lem}

\begin{proof}
Since $\beta_i\to 0$ and $\|\Inj_{i,i+1}\|\le 1$, we have $W_1P_N-W_1\to 0$, with $P_N$ the natural projection onto the first $N$ coordinates. The ideal of superstrictly singular operators being norm closed, it is enough to show that each $W_1P_N$ is superstrictly singular. But $W_1P_N$ is obtained by adding as direct summands zero operators to the operator
\[
W'_N=\beta_1\Inj_{p_1,p_{2}}\oplus \cdots \oplus \beta_N\Inj_{p_N,p_{N+1}}: \bigoplus_{i=1}^N J_{p_i} \to \bigoplus_{i=2}^{N+1} J_{p_i}.
\]
Since by Theorem~\ref{th:main} all summands are superstrictly singular, it follows that $W'_N$ is superstrictly singular. Therefore $W_1$ is superstrictly singular.
\end{proof}

The arguments in~\cite{R} yield then the next theorem.
\begin{Thm}\label{th:read}
Read's operator $T$ is superstrictly singular.
\end{Thm}
Thus $T$ is a strictly singular operator without nontrivial invariant subspaces.

\section{The technical lemma}\label{se:proof}

Fix a natural number $N$, and denote $\BB_N=\{x\in \RR^N: \|x\|_\infty\le1\}$,
$\Sigma^{N-1}=\{x\in \RR^N: \|x\|_\infty=1\}$.
For $k\ge 1$ we define
\begin{align*}
\Gamma(k)&=\{x\in \BB_{N}: x\mbox{ has at least $k$ alternating coordinates $\pm1$}\},\\
A_+(k)&=\{x\in \BB_{N}: x\mbox{ has at least $k$ alternating coordinates $\pm1$, starting with 1}\},\\
A_-(k)&=-A_+(k).
\end{align*}
Put also $A_+(0)=A_-(0)=\Gamma(0)=\BB_N$.
For $k\ge1$, $\Gamma(k), A_{\pm}(k)\subset \Sigma^{N-1}$ and
we have 
\begin{align*}
	A_+(k)\cup A_-(k)&= \Gamma(k),\\%\label{eq:0}\\ 
 A_+(k)\cap A_-(k)&= \Gamma(k+1).
	%\label{eq:1}
\end{align*}
Note that the first relation above is true also for $k=0$.

We start with a simple lemma.

\begin{Lem}\label{le:0}
Suppose $p$ is a real polynomial of degree $m$, and there are $m+2$ real numbers $t_1<t_2<\dots<t_{m+2}$, such that $p(t_i)\ge0$ for $i$ odd and $p(t_i)\le 0$ for $i$ even. Then $p\equiv 0$.
\end{Lem}

\begin{proof}
We do induction with respect to $m$. If $m=0$, the result is obvious. If the lemma has been proved up to $m-1$, and  $p$ is a polynomial of degree $m$, then $p$ has at least one real root $s$. We write $p(t)=(t-s)q(t)$, and $q$ (or $-q$) has a similar property, with respect to at least $m-1$ values $t_i$---so we can apply induction.
\end{proof}

\begin{Lem}\label{le:1}
There exists a sequence of subspaces $\pi_k\subset \RR^N$, $\pi_k\supset \pi_{k+1}$, $\dim\pi_k=N-k$, such that, if $P_k$ is the orthogonal projection onto $\pi_k$, then $P_k|A_+(k)$ is injective.
\end{Lem}

\begin{proof} 
For $1\le j\le N$ we define the vectors $\zeta^j\in \RR^N$ by the formula $\zeta^j_i=i^{j-1}$. One checks easily that the $\zeta^j$'s are linearly independent. Define $\pi_0=\RR^N$, and, for $k\ge 1$, $\pi_k=(\Span\{ \zeta^1,\dots, \zeta^k\})^\perp$.

Suppose that $x,y\in A_+(k)$, and $P_k(x)=P_k(y)$. There exist scalars $\alpha_1,\dots, \alpha_k$, such that $x-y=\sum_{j=1}^k\alpha_j\zeta^j$. 
We have indices $1\le r_1<\dots<r_k\le N$ and $1\le s_1<\dots<s_k\le N$, such that $x_{r_l}=y_{s_l}=(-1)^{l-1}$. It follows that $x_{r_l}-y_{r_l}\ge 0$ for $l$ odd and $\le0$ for $l$ even, while $x_{s_l}-y_{s_l}\le 0$ for $l$ odd and $\ge0$ for $l$ even.

Let the polynomial $p$ of degree $k-1$ be given by $p(t)=\sum_{j=1}^k\alpha_j t^{j-1}$.
If $r_l=s_l$ for all $l$, we obtain \[\sum_j \alpha_j\zeta^j_{r_l}=\sum_j\alpha_j {r_l}^{j-1}=0\]
for all $l=1,\dots k$. Thus $p$  has $k$ distinct zeros; it must be identically 0, whence $x=y$.

Suppose now that we have $r_l\not=s_l$ for at least one index $l$. We claim then that among the union of the indices $r_l$ and $s_l$ we can find $\iota_1<\iota_2<\dots<\iota_{k+1}$, such that $x_{\iota_l}-y_{\iota_l}$ have alternating signs. This can be achieved by induction with respect to $k$. For $k=1$ we must have $r_1\not= s_1$, so we may take $\iota_1=\min\{r_1, s_1\}$, $\iota_2=\max\{r_1, s_1\}$. 
For  $k>1$, there are two cases. If $r_1=s_1$, we take $\iota_1=r_1=s_1$ and apply the induction hypothesis to obtain the rest. If $r_1\not=s_1$, we take $\iota_1$ whichever is the first among them, $\iota_2$ as the other one, and then we continue ``accordingly" to $\iota_2$ (that is, taking as $\iota$'s the rest of $r$'s if $\iota_2=r_1$ and the rest of $s$'s if $\iota_2=s_1$).

Now, the way  $\iota_l$ have been chosen implies that $p(t)$ defined above satisfies the hypotheses of Lemma~\ref{le:0}: it has degree $k-1$ and the values it takes in $\iota_1,\dots,\iota_{k+1}$ have alternating signs. It must then be identically 0, which implies $x=y$. 
\end{proof}

Since $A_-(k)=-A_+(k)$, it follows that $P_k|A_-(k)$ is also injective.

\begin{Lem}\label{le:rec}
If $\pi_k, P_k$ are obtained in Lemma~\ref{le:1}, then
\[
\Delta_k:=P_k(\Gamma(k))%=P_k(A_-(k)) 
\]
is a balanced, convex subset of $\pi_k$, with $0$ as an interior point (in $\pi_k$). Moreover, $\Delta_k=P_k(A_-(k)=P_k(A_-(k))$ and $\partial\Delta_k=P_k(\Gamma(k+1))$ (the boundary in the relative topology of $\pi_k$).
\end{Lem}

\begin{proof}
We will use induction with respect to $k$. The statement is immediately checked for $k=0$ (note that $P_0=I_{\RR^N}$ and $\partial \Delta_0=\Sigma^{N-1}=\Gamma(1)$).

%%%%%%%%%

Assume the statement true for $k$; we will prove its validity for $k+1$.
By the induction hypothesis, we have
\[
\Delta_{k+1}=P_{k+1}P_k(\Gamma(k+1))=P_{k+1}\partial\Delta_k= P_{k+1}\Delta_k
\]
and is therefore a balanced, convex subset of $\pi_{k+1}$, with 0 as an interior point.

Take then $y\in \Int{\Delta}_{k+1}$. Suppose  $P_{k+1}^{-1}(y)\cap\partial\Delta_k$ contains a single point. Then
$P_{k+1}^{-1}(y)\cap\Delta_k$ also contains a single point, and therefore $P_{k+1}^{-1}(y)\cap \pi_k$ is a  support line for the convex set $\Delta_k$. This line is contained in a support hyperplane (in $\pi_k$); but then the whole of $\Delta_k$ projects onto $\pi_{k+1}$ on one side of this hyperplane, and thus $y$ belongs to the boundary of this projection. Therefore $y$ cannot be in $\Int{\Delta}_{k+1}$.

The contradiction obtained shows that $P_{k+1}^{-1}(y)\cap\partial\Delta_k$ contains at least two points. But 
\[
\partial\Delta_k=P_k(\Gamma(k+1))= P_k(A_+(k+1))\cup P_k(A_-(k+1))
\]
whence
\[
P_{k+1}(\partial\Delta_k) = P_{k+1}(A_+(k+1))\cup P_{k+1}(A_-(k+1)).
\]
Since $P_{k+1}$ restricted to each of the two terms in the right hand side is injective by Lemma~\ref{le:1}, it follows that $P_{k+1}^{-1}(y)\cap\partial\Delta_k$ is formed by exactly two points, one in $A_+(k+1)$ and the other in $A_-(k+1)$.
In particular, $\Int{\Delta}_{k+1}\subset P_{k+1}(A_\pm(k+1))$. But, $\Delta_{k+1}$ being a closed convex set with nonempty interior, it is the closure of its interior $\Int{\Delta}_{k+1}$; since the two sets on the right are closed, we have actually $\Delta_{k+1}= P_{k+1}(A_\pm(k+1))$.

We want to show now that  $\partial\Delta_{k+1}=P_{k+1}(\Gamma_{k+2})$.
Suppose first that $y\in P_{k+1}(\Gamma_{k+2})= P_{k+1}(A_+(k+1)\cap A_-(k+1))$. Then the line $\ell\subset\pi_k$ orthogonal in $y$ to $\pi_{k+1}$ cannot have other points of intersection with $A_+(k+1)$ or with $A_-(k+1)$, since $P_{k+1}$ is injective on these two sets. Therefore $\ell$ is a line of support of $\Delta_{k+1}$, and is contained in a support hyperplane (inside $\pi_k$), whence $y\in\partial\Delta_{k+1}$.

Conversely, take $y\in\partial\Delta_{k+1}=\partial(P_{k+1}(\Delta_k)$.  Take $z_+\in A_+(k+1)$, $z_-\in A_-(k+1)$, such that $P_{k+1}(z_+)=P_{k+1}(z_-)=y$. We have then $P_{k}(z_+)\in\partial\Delta_k$ (if $P_{k}(z_+)\in \Int{\Delta}_k$, then $P_{k+1}(z_+)=P_{k+1}(P_{k}(z_+))$ must be in the interior of ${P_{k+1}\Delta_k}$, which is $\Int{\Delta}_{k+1}$). Similarly, $P_{k}(z_-)\in\partial\Delta_k$. 

If $P_{k}(z_+)\not= P_{k}(z_-)$, then $P_{k+1}$ applied to the whole segment $[P_{k}(z_+), P_{k}(z_-)]$ is equal to $y$. Therefore the segment belongs to $\partial\Delta_k$. Since $\partial\Delta_k= P_k(A_+(k+1)\cup A_-(k+1))$, there exist two values $x_1, x_2$ either both in $A_+(k+1)$ or both in $A_-(k+1)$, such that $P_{k}x_1, P_{k}x_2\in [P_{k}(z_+), P_{k}(z_-)]$, and thus $P_{k+1}x_1=P_{k+1}x_2=y$. This contradicts the injectivity of $P_{k+1}$ on  $A_+(k+1)$.

Therefore $P_{k}(z_+)= P_{k}(z_-)$. But $z_+$ and $z_-$ both belong to $A_+(k)$, on which $P_k$ is injective. It follows that $z_+=z_-\in A_+(k+1)\cap A_-(k+1)=\Gamma(k+2)$, and $P_{k+1}z_+=y$. This ends the proof.
\end{proof}

The main consequence of Lemma~\ref{le:rec}, in combination with Lemma~\ref{le:1}, is the fact that 
the linear map $P_{k-1}$ maps  homeomorphically $\Gamma(k)$ into $\partial\Delta_{k-1}$, which is the boundary of a convex, balanced set, containing 0 in its interior.

\begin{Cor}\label{co:1}
If $X\subset \RR^N$ is a subspace of dimension $k$, then $X\cap \Gamma(k)\not=\emptyset$.
\end{Cor}

\begin{proof}
As noted above, $P_{k-1}$ maps homeomorphically $\Gamma(k)$ onto the boundary of a convex, balanced set, containing 0 in its interior. Composing it with the map $x\mapsto \frac{x}{\|x\|}$, we obtain a homeomorphic map $\phi$ from $\Gamma(k)$ to $S^{N-k}$, which satisfies the relation $\phi(-x)=-\phi(x)$.

If $X\cap \Gamma(k)=\emptyset$, then the projection of $\Gamma_k$ onto $X^\perp$ does not contain 0. Composing this projection with the map $x\mapsto \frac{x}{\|x\|}$, we obtain a continuous map from $\psi:\Gamma(k)\to S^{N-k-1}$, that satisfies $\psi(-x)=-\psi(x)$. Then the map $\Phi:=\psi\circ\phi^{-1}:S^{N-k}\to S^{N-k-1}$ is continuous and satisfies $\Phi(-x)=-\Phi(x)$. This is however impossible: it is known that such a map does not exist (see, for instance,~\cite{F}). Therefore we must have $X\cap \Gamma(k)\not=\emptyset$.
\end{proof}

Finally, we extend this result to the infinite dimensional Banach space $c_0$ of real-valued sequences that tend to zero, as stated in Lemma~\ref{le:pm1}.

%\begin{Thm}\label{le:pm1}
%If $X\subset c_0$ is a subspace of dimension $k$, then there is $x\in X$, $\|x\|=1$ and indices $i_1 <\dots <i_k$, such that $x_{i_j}=(-1)^j$.
%\end{Thm}

\begin{proof}[Proof of Lemma~\ref{le:pm1}]
Denote by $P_N$ the orthogonal projection onto the first $N$ coordinates. Since $I-P_N\to 0$ on each element of $c_0$, and $X$ is finite dimensional, a standard argument shows that $I-P_N|X:X\to c_0$ tends to 0 uniformly. In particular, $P_N|X$ is bounded below for large $N$, and thus
$X_N=P_NX$ has eventually dimension $k$.

Applying Corollary~\ref{co:1} to $X_N$ yields vectors $y^N=P_Nx^N$, $x^N\in X$, $\|y_N\|=1$, whose coordinates have $k$ alternating $1$ and $-1$. We have $x^N=y^N+(I-P_N)x^N$, whence $\|x^N\|\le 1+\|(1-P_N)x^N\|$; therefore, if we take $\|I-P_N\|<1/2$, we obtain that the sequence $x^N$ is bounded. Take $x\in X$ to be the limit of a convergent subsequence.
Since 
\[
\|(I-P_N)x^N\|\le \|(I-P_N)x\|+\|I-P_N\|\cdot \|x-x^N\|,
\]
the indices of the coordinates of $y^N$ of modulus 1 are bounded by an absolute constant. By passing to a subsequence we may assume that they are constant. It follows then that $x$ has the corresponding $\pm1$ on those coordinates. 
\end{proof}

\end{document}